\documentclass[aps,pra,amssymb,amsmath,showkeys]{revtex4-1}

\usepackage{bm}
\usepackage{graphicx}

\newtheorem{thm}{Theorem}
\newtheorem{lmt}[thm]{Lemma}
\newtheorem{cor}[thm]{Corollary}

\newtheorem{rmk}{Remark}

\newtheorem{prf}{Proof}

\newcommand*{\qed}{\hbox{}\hfill$\Box$}

\begin{document}

\title{Spectrum of some arrow-bordered circulant matrix}

\author{Wojciech Florek}
\affiliation{Adam Mickiewicz University, Faculty of Physics,
ul.\ Uniwersytetu Pozna\'nskiego~2, 61-614 Pozna\'n, Poland}
\email{wojciech.florek@amu.edu.pl}

\author{Adam Marlewski}
\affiliation{Pozna\'n University of Technology, Institute of Mathematics,
ul.\ Piotrowo 3A, 60-965 Pozna\'n, Poland}

\date{\today}

\begin{abstract}
Given a circulant matrix $\mathrm{circ}(c,a,0,0,...,0,a)$, $a\ne 0$, of 
order~$n$, we ``border'' it from left and from above by constant column and 
row, respectively, and we set the left top entry to be $-nc$. This way we get 
a~particular title object, an example of what we call an \textit{abc matrix\/}, 
or an \textit{arrow-bordered circulant (matrix)\/}. We find its eigenpairs and 
we discuss its spectrum with stress on extreme eigenvalues and their bounds. 
At last we notice its relation to a~weighted wheel graph.  
\end{abstract}


\maketitle


\section{Introduction and motivation}\label{intro}
By ``gluing'' a~circulant matrix and an arrow(head) matrix we get an 
\textit{abc matrix\/} (``abc'' acronymizes ``arrow-bordered circulant''). We 
notice that an abc matrix is the adjacency matrix of a~wheel graph whose both 
vertices and edges have weights. By elementary methods we find eigenpairs of 
a~regular abc matrix and we analyse its spectrum. The work with matrices we 
consider is motivated by problems such as that treated by \citet{Schm03}, 
where elements of spectral graph theory are applied to determine the lowest 
energy configuration(s) of finite classical Heisenberg spin systems. There 
exists a~weighted adjacency matrix of a~(simple) graph representing such system 
of interacting localized spin vectors (and this graph, $W_n$, is discussed in 
Section~\ref{wheel} below), weights are assigned to both vertices and edges; in 
addition, this matrix is dressed with a~``gauge'' vector spread over the main 
diagonal (the notion ``dressing with a~vector'' is explained in 
Section~\ref{cncl}). The crucial point in the search of desired configurations 
is to find the minimum eigenvalue of the matrix at hand---the sought solution 
is determined by the maximum, determined with respect to the gauge, of the 
minimum eigenvalue. The other eigenvalues have no clear physical meaning, 
nevertheless properties of the whole spectrum are discussed below.

There is an extensive literature on spectral graph theory. Numerous relations 
between graphs and their 0-1 adjacency matrices are presented, among others, in 
monographs \citet{Bigg74}, \citet{Bond76}, \citet{Cvet88}, \citet{Brou12}, and 
\citet{Nica18}, where, in particular, the Cauchy theorem on interlacing 
eigenvalues of a~symmetric matrix and its submatrix, and Weyl upper bound for 
the eigenvalues of a~sum of two symmetric matrices are recalled. \citet{Zhan09} 
discussed adjacency matrices of graphs obtained by attaching some pendent edges 
to complete graphs ($K_n$) and to cyclic graphs ($C_n$), whereas \citet{DasK13} 
dealt with eigenvalues of a~friendship graph ($F_n$). Introduced in this work 
an abc matrix $M_n$ (see Section~\ref{abcm}) is the adjacency matrix of 
a~regularly weighted wheel graph, the notion we define in Section~\ref{wheel}. 
Graphs and their adjacency matrices are widely used among others in chemistry, 
physics and telecommunications \citep[a~survey is presented in][]{Cvet11}, 
and it mostly concerns unweighted graphs. Examples where edge- or 
vertex-weighted graphs are dealt with include \citet{Klav97}, \citet{Lato03}, 
\citet{Ashr05}, \citet{Mirc15}, \citet{Jooy16}, \citet{CaiH17}, and 
\citet{Anto18}.

The paper is organized as follows: In Section~\ref{base}, the most important 
properties of circulant and arrowhead matrices are recalled. The title object 
of this work is presented in Section~\ref{abcm}, whereas its eigenpairs are 
determined in Section~\ref{egabc}. Main theorems are included in three parts of 
Section~\ref{analyz}. In Section~\ref{wheel} it is shown hat the abc matrices 
discussed in this paper can be treated as weighted adjacency matrices of wheel 
graphs. Finally, Section~\ref{cncl} is devoted to some final remarks.

\section{Basics on circulant matrices and arrow(head) matrices}\label{base}
In linear algebra a~matrix $[a_{j,k}]$ is called 
\textit{Toeplitz},\footnote{In 1911, in the paper \textit{\"{U}ber allgemeine 
lineare Mittelbildungen\/} (published in annals \textsl{Prace 
Matematyczno-Fizyczne}; they appeared in Warsaw, its 48 volumes were published 
in years 1888--1952), Otto Toeplitz (1881--1940), when discussing the 
summability of series, treated specific matrices, in particular triangular 
ones. Two years later his result was generalized (and it is now called the 
Silverman-Toeplitz theorem), and infinite matrices satisfying some conditions 
became to be referred to as Toeplitz; note that these Toeplitz matrices are 
objects distinct from that considered in linear algebra.} or 
\textit{constant-diagonal}, if its each descending diagonal from left to right 
is constant, so if $a_{j,k} = t_{j-k}$ (and we say that quantities $t_m$'s 
determine it). Such infinity matrices, $T = [t_{j-k}]_{j,k = 1,2,3,\ldots}$, 
were first investigated by Otto Toeplitz, he showed that the matrix~$T$ defines 
a~bounded linear operator, which transforms a~sequence $x\in \ell^2$ in 
$Tx = y$, where $y_k =\sum_{k=0,1,2,\ldots} t_{j-k} x_j$, iff $t_m$ are Fourier 
coefficients of an appropriate function. In fact, Toeplitz proved it in the 
special case of symmetric matrix ($t_{-m} = t_m$), and several decades later 
his result was extended to the general case by \citet{Hart54}. Below we deal 
with finite Toeplitz matrices. By definition, 
a~$(m,n)$-\textit{Toeplitz matrix\/} has $m$~rows and $n$~columns, its 
top/first row is $[t_0,t_1,t_2,\ldots,t_{n-2},t_{n-1}]$, the second row is 
$[t_{-1},t_0,t_1,\ldots,t_{n-3},t_{n-2}]$, etc. In particular, a~square 
Toeplitz matrix, $[a_{j,k}]_{j,k=1,\ldots,n}$, has its last/$n$th row 
$[t_{-(n-1)},t_{-(n-2)},t_{-(n-3)},\ldots,t_{-1},t_0]$. So, any Toeplitz matrix 
of order~$n$ (and of $2n-1$ degrees of freedom) is fully determined by the 
vector 
$\bm{t}=[t_{-(n-1)},t_{-(n-2)},t_{-(n-1)},\ldots,t_{-1},t_0,t_1,\ldots,t_{n-1}]$, 
and is denoted by $\mathrm{toep}(\bm{t})$.

For $a_{j,k}=t_{|j-k|}$ we have a~\textit{symmetric Toeplitz matrix\/} of 
order~$n$, and we denote it by the same symbol, $\mathrm{toep}(\bm{t})$, where 
$\bm{t}=[t_0,t_1,\ldots,t_{n-1}]$ is its first row. A~review of results on 
eigenvalues of symmetric Toeplitz matrices is given in \citet{Dels84}, 
\citet{Reic92}, \citet{Melm01}; in \citet{Andr73}, \citet{Cant76}, 
\citet{LiuZ03}, \citet{AbuJ04}, \citet{Kato14}, and \citet{Brua17} there are 
mostly discussed, general or particular, centrosymmetric 
matrices.\footnote{A~matrix $a=[a_{j,k}]_{j,k=1,2,\ldots,n}$ is 
said to be \textit{centrosymmetric}, or \textit{cross-symmetric}, if it is 
symmetric about its center, $a_{j,k}=a_{n-1+j,n-1+k}$. A~necessary and 
sufficient condition a~matrix~$[a_{jk}]$ to be centrosymmetric is its 
commutation with so-called exchange matrix~$J$ (by definition, all entries 
of~$J$ are~0, but $J_{j,n+1-j}=1$ for all $j=1,2,\ldots,n)$. Thus every 
symmetric Toeplitz matrix is centrosymmetric. The centrosymmetricity is one of 
patterns  of symmetry \citep[see, e.g.,][]{Weav85,Pres98,Tren04}.}

A~square Toeplitz matrix whose $j$th row is the cyclic shift of the top/first 
row, $\bm{c}=[c_0,c_1,\ldots,c_{n-1}]$ by $j-1$ positions to the right, is 
called \textit{circulant (matrix)\/} and denoted by $\mathrm{circ}(\bm{c})$, or 
$\mathrm{circ}(c_0,c_1,\ldots,c_{n-1})$. First circulant matrices were studied 
by Eug\`{e}ne Catalan in \textit{Recherches sur les d\'{e}terminants\/} (1846), 
William Spottiswoode in \textit{Elementary theorems relating to determinants\/} 
(1856) and Alphonse Legoux in \textit{Application d'un d\'{e}terminant\/} 
(1883). Circulant matrices appear in numerous problems 
\citep[see, e.g.,][]{Gris06,Diac12,Olso14}. A~classical book on circulant 
matrices is \citet{Davi94}, whereas more recent texts dedicated to them and 
their generalizations are \citet{Gray06,Fuhr12,KraI12,Bose12}. They all cite 
the following basic result on spectral properties of arbitrary circular matrix 
\begin{thm}\label{basicthm}
Given $\bm{c}:=[c_0,c_1,\ldots,c_{n-1}]$, eigenpairs of $\mathrm{circ}(\bm{c})$
are $(\lambda_k,\bm{v}_k)$, $k=0,1,2,\ldots,n-1$, where
\begin{eqnarray*}
 \lambda_k &:=& c_0+c_1\omega_n^k+c_2\omega_n^{2k}
    +\cdots+c_{n-1}\omega_n^{(n-1)k},\\
 \bm{v}_k&:=& [1,\omega_n^k,\omega_n^{2k},\omega_n^{3k},
   \ldots,\omega_n^{(n-1)k}]^\mathsf{T},\\
 \omega_n &:=& \exp(2\pi\mathrm{i}/n);
\end{eqnarray*}
$\omega_n$ is the first prime root of degree~$n$ of the unity and, as usual, 
the superscript ``\textsf{T}'' stands for the transpose.
\end{thm}

A~square matrix containing zeros in all its entries except for the first row, 
first column, and main diagonal, is called (because of the pattern its nonzero 
entries form) an \textit{arrowhead}, or (\textit{left up}) \textit{arrow}, 
\textit{matrix} \citep[e.g.][]{Shen09}. The notion ``arrow matrix'' embraces 
also all matrices similar to any arrowhead matrix via a~symmetric permutation 
\citep[but there are also in use more refined names, e.g., down-arrow matrix, 
see][]{Parl09}. 

Diagonal matrices are arrow matrices, so examples of arrow matrices appeared 
when, in 1850 in \textit{Additions to the articles, ``On a New Class of 
Theorems'', and ``On Pascal's Theorem''}, Joseph James Sylvester coined the 
term ``a~matrix''. A~good compendium on arrow matrices is \citet{OLea90}. 
Examples of centrosymmetrical arrow matrices include matrices 
$A_n:=A_n(h,b,d)\in\mathbb{R}^{(n+1)\times(n+1)}$, $n>0$, where $h$~is the most 
left top entry (one can refer it to as a~headpoint entry), all other bordering 
entries (i.e., elements in the top row and in the first column) are~$b$, and 
all other diagonal elements are~$d$. So, it has a~form as follows:
\begin{equation}\label{ram}
A_n(h,b,d)=\left[\begin{array}{c|c} 
   h & \bm{b}^\mathsf{T} \\ \hline \bm{b} & D
\end{array}\right],
\end{equation}
where $h\in\mathbb{R}$, $\bm{b}=[b,b,\ldots,b]^\mathsf{T}\in\mathbb{R}^n$ is 
a~constant vector, $D=\mathrm{diag}(d,\ldots,d)$ is a~diagonal matrix with 
$d\in\mathbb{R}$. Such arrow matrices are referred to as \textit{regular\/} 
ones. Their eigenvalues are given by
\begin{thm}\label{ahme}
Let $A_n\in\mathbb{R}^{(n+1)\times(n+1)}$ be a~regular arrowhead matrix of the 
form~(\ref{ram}). Then, its eigenvalues are
\begin{subequations}\label{Anev}
\begin{eqnarray}
 \lambda_\pm &=\quad \left(h+d\pm\sqrt{\Delta_n}\right)/2, 
 \qquad  & \mathrm{for}\quad n\ge 1,\label{Anevpm}\\
 \lambda_k &=\quad d, \qquad \qquad\qquad\qquad\;
    & \mathrm{for}\quad n>1\quad \mathrm{and}\quad 
    k=1,2,\ldots,n-1,\label{Anevk} 
\end{eqnarray}
\end{subequations}  
where $\Delta_n=(h-d)^2+4nb^2$.
\end{thm}
\begin{prf} 
Let $\alpha_n(\lambda)$ denote the characteristic polynomial of the 
matrix~$A_n$. Then 
\begin{equation}\label{alph1}
\alpha_1(\lambda)=(h-\lambda)(d-\lambda)-b^2
  =\lambda^2-(h+d)\lambda+hd-b^2.
\end{equation}
For $n>1$ we have
\begin{equation}\label{alphn}
  \alpha_n(\lambda)=(d-\lambda)^{n-1}\left((h-\lambda)(d-\lambda)-nb^2\right).
\end{equation}
Indeed, applying the Laplace expansion we have inductively
\[
\alpha_{n+1}(\lambda)= (d-\lambda)\alpha_n(\lambda)-(d-\lambda)^n b^2
 =(d-\lambda)^n\left((h-\lambda)(d-\lambda)-(n+1)b^2\right).
\]
The first factor in Eq.~(\ref{alphn}) proves (\ref{Anevk}), whereas the second 
part and Eq.~(\ref{alph1}) lead to a~quadratic equation with the discriminant 
$\Delta_n=(h-d)^2+4nb^2$, $n>0$, and, therefore, the other two eigenvalues 
are given by Eq.~(\ref{Anevpm}), which completes the proof.\qed
\end{prf}

For $n>1$ this theorem can be proved applying Corollary~4 stated by 
\citet{Shen09}. Note that $\lambda_-<\lambda_+$ if $b\ne 0$, despite values of 
$h$~and~$d$. Immediately from Theorem~\ref{ahme} it follows
\begin{cor}\label{specA}
For any real numbers~$h$, $d$ and $b\ne 0$ the spectrum 
$\sigma\bigl(A_n(h,b,d)\bigr)$ is of cardinality 
\begin{equation}\label{specAe}
\left|\sigma\bigl(A_n(h,b,d)\bigr)\right|= 
\begin{cases} 2, & \mathrm{for}\qquad n=1,\\
              3, & \mathrm{otherwise}.
\end{cases}
\end{equation}
\end{cor}

\section{Arrow-bordered circulant matrix}\label{abcm}
Let us define an \textit{arrow-bordered\/} (or \textit{arrowly bordered}) 
\textit{circulant matrix}, or an \textit{abc matrix} for short, as 
a~circulant matrix expanded on its left with a~column, and on its top with 
a~row. So an abc matrix $m_n=[m_{j,k}]_{j,k=0,1,...,n}$, $n>0$, is block 
structured as follows:
\begin{equation}\label{abcdef}
  m_n = \left[\begin{array}{r|r} 
    h & \bm{r}^\mathsf{T} \\ \hline \bm{b} & t_n 
  \end{array}\right],
\end{equation}
where a~scalar~$h$ is proposed to be referred to as a~headpoint number (or 
a~tip number), $\bm{b}$~and~$\bm{r}^\mathsf{T}$ are (bordering) vectors, and 
$t_n$~is a~circulant matrix of order~$n$. Symmetric arrow-bordered matrices 
(i.e., for $\bm{b}=\bm{r}$) with 
a~diagonal~$t_n=\mathrm{diag}(t_1,t_2,\ldots,t_n)$, called \textit{headarrow 
matrices}, are treated in \citet{OLea90}, \citet{Pick07}, \citet{Shen09}, and 
\citet{Jako15}. Below we deal with more general case, namely with 
$t_n=\mathrm{circ}(c,a,0,\ldots,0, a)$, where $a\ne0$ in a~general case. 

A~symmetrical abc matrix with a~constant vector 
$\bm{b}=[b,b,\ldots,b]^\mathsf{T}=\bm{r}$ is referred to as 
a~\textit{regularly arrow-bordered matrix}, or a~\textit{regular abc matrix}. 
In the next we consider \textit{traceless\/} regular abc matrices, so matrices 
$M_n:=m_n(a,b,c)$, where $a$, $b$, $c$ are real numbers, 
$t_n=\mathrm{circ}(c,a,0,0,\ldots,0,0,a)$, and $h=-nc$ (this choice makes that 
$M_n$ is traceless). For example (zero entries are marked by dots),
\[
M_6 = m_6(a,b,c) =
 \left[\begin{array}{r|r}-6c& \bm{b}^\mathsf{T}\\ \hline \bm{b} & t_6
 \end{array}\right] 
= \left[\begin{array}{r|*{6}{r}}
   -6c & b & b & b & b & b & b \\ \hline
     b & c & a & \cdot & \cdot & \cdot & a \\ 
     b & a & c & a & \cdot & \cdot & \cdot \\
     b & \cdot & a & c & a & \cdot & \cdot \\ 
     b & \cdot & \cdot & a & c & a & \cdot \\
     b & \cdot & \cdot & \cdot & a & c & a \\ 
     b & a & \cdot & \cdot & \cdot & a & c 
\end{array}\right].
\]

\begin{rmk}\label{trabc}
Some obvious but important remarks are in place.
\begin{enumerate}
\item \label{rmka}
  With $b=0$ we have a~trivial situation, $m_n(a,0,c)=\mathrm{diag}(-nc,t_n)$, 
  so in this case the spectrum 
\[\sigma\bigl(m_n(a,0,c)\bigr)=\{-nc\}\cup\sigma(t_n).\]
  Therefore in the next we deal with $b\ne 0$.
\item \label{rmkb}
  Since $m_n(a,b,c) = b\,m_n(a/b,1,c/b)$, we can deal with abc matrices with 
  fixed $b=1$ or $b=-1$. Obviously, this does not restrict the generality of 
  considerations, and we mainly discuss these cases, i.e., $b=\pm1$.
\item \label{rmkc}
  The special case $a=0$ gives $m_n(0,b,c)=A_n(-nc,b,c)$ and we treated it in 
  Theorem~\ref{ahme}.
\end{enumerate}
\end{rmk}

\section{Eigenpairs of an abc matrix}\label{egabc}
All matrices that we treat in the next are regular abc matrices, and we pay our 
attention to traceless ones. Nevertheless, for reasons that appear clear later, 
we need to treat two kinds of abc matrices of order two and three ($n=1,2$, 
respectively), namely that defined by~(\ref{m12tilde}) and that defined 
by~(\ref{m12notilde}); as we will see, formulas (\ref{m12tilde}) perfectly 
match the general case ($n>2$, $a\ne0$), while formulas (\ref{m12notilde}) do 
not (but, surprisingly, go well with Theorem~\ref{egpair}). In corresponding 
definitions we take different vectors~$\bm{t}$.

We can take $\bm{t}=[c]$ and $\bm{t}=[c,a]$ and associate to them matrices 
$\tilde{t}_1=\mathrm{circ}(c)$ and $\tilde{t}_2=\mathrm{circ}(c,a)$, 
respectively. This way we get
\begin{equation}\label{m12tilde}
\widetilde{M}_1=
\left[ \begin{array}{rr}
      -c & b \\ b & c 
\end{array}\right]=A_1(-c,b,c)
\qquad\mathrm{and}\qquad
\widetilde{M}_2 = \left[ \begin{array}{rrr}
      -2c & b & b\\ b & c & a \\ b & a &c
\end{array}\right].
\end{equation}
In the above definitions tilded symbols are used, since these matrices do not 
observe general formulas given below, in Theorem~\ref{egpair}. Note, that in 
this approach the matrix~$\widetilde{M}_1$ does not depend on the 
parameter~$a$. It is easy to see
\begin{thm}[Spectra of $\widetilde{M}_1$ and $\widetilde{M}_2$] 
\mbox{~}\label{thmtilde}
\begin{enumerate}
  \item $\sigma\bigl(\widetilde{M}_1\bigr)=\left\{\pm\sqrt{b^2+c^2}\right\}$, 
    so  $\left|\sigma\bigl(\widetilde{M}_1\bigr)\right|=2$. 
  \item $\sigma\bigl(\widetilde{M}_2\bigr)=\{\lambda_-,\lambda_+,\lambda_1\}$,
    where
\begin{eqnarray*}
  \lambda_\pm &:=& \frac{1}{2}\left(a-c\pm\sqrt{\tilde{\Delta}}\right) \\
  \lambda_1   &:=& c-a,
\end{eqnarray*}
and $\tilde{\Delta}=(a+3c)^2+8b^2$. Therefore,
\[
 \left|\sigma\bigl(\widetilde{M}_2\bigr)\right|=\begin{cases}
        2, & \mathrm{if}\quad a\ne0\quad\mathrm{and}\quad c=(a^2-b^2)/(3a),\\
        3, & \mathrm{otherwise}.
  \end{cases}
\]
\end{enumerate}
\end{thm}
Let us pay attention that
\begin{enumerate}
  \item For $c=(a^2-b^2)/(3a)$, $a\ne0$, there is 
    $\sigma\bigl(\widetilde{M}_2\bigr)=\{s,-s/2\}$, where $s:=2(2a^2+b^2)/(3a)$ 
    is the~single eigenvalue. 
  \item The lines $\lambda=c+a$ and $\lambda=-2c$ are asymptotes to both curves 
    $\lambda=\lambda_-$ and $\lambda=\lambda_+$.
  \item Since $\tilde{\Delta}>0$ for $b\ne0$, then $\lambda_-<\lambda_+$.
\end{enumerate}

The other possibility is to demand that a~sum of entries of the vector~$\bm{t}$ 
equals $c+2a$. Within this approach $t_1:=\mathrm{circ}(c+2a)$ and 
$t_2:=\mathrm{circ}(c,2a)$, so
\begin{equation}\label{m12notilde}
M_1=
\left[ \begin{array}{rr}
      -c & b \\ b & c+2a 
\end{array}\right]=A_1(-c,b,c+2a)
\qquad\mathrm{and}\qquad
M_2 = \left[ \begin{array}{rrr}
      -2c & b & b\\ b & c & 2a \\ b & 2a &c
\end{array}\right].
\end{equation}
The spectrum and its properties of the matrix~$M_1$ can be determined from 
Theorem~\ref{ahme}. For the matrix~$M_2$ Theorem~\ref{thmtilde} may be applied 
replacing the parameter~$a$ by~$2a$. The notions introduced below in 
Eqs.~(\ref{mndef}) can be formally applied to the matrix~$M_1$, but this matrix 
has a~nonvanishing trace and is omitted in the further discussion in 
Sections~\ref{analyz} and~\ref{wheel}. On the other hand, the matrix~$M_2$ 
defined above obeys all assumptions and, therefore, can be included in general 
considerations, so it is assumed hereafter that the matrix $m_2(a,b,c)$ is 
constructed with $t_2=\mathrm{circ}(c,2a)$. In Section~\ref{wheel}, where 
a~relation of abc matrices to wheel graphs is presented, both matrices, 
$\widetilde{M}_2$ and $M_2$, are taken into account.

Completed the discussion on $n\in\{1,2\}$, we fix a~natural number $n>2$, take 
arbitrary real numbers $a$, $b\ne0$, $c$, and denote
\begin{subequations}\label{mndef}
\begin{eqnarray}
  M_n &:=& m_n(a,b,c),\label{mnM}\\
 \varphi_n &:=& 2\pi/n,\qquad\mathrm{so}\quad 
    \omega_n = \exp(\mathrm{i}\varphi_n), \\
 \Delta_n:= \Delta(a,b,c,n) &:=&  \bigl(2a+(n+1)c\bigr)^2 + 4nb^2,\label{mnD}\\
 \beta_{n;\pm}:= \beta_\pm(a,b,c,n) &:=& 
    -\Bigl(2a+(n+1)c\mp\sqrt{\Delta_n}\Bigr)/(2b), \label{mnB}  \\ 
 \lambda_{n;\pm} := \lambda_\pm(a,b,c,n) &:=& 
    b\beta_{n;\pm} +2a+c
   = \Bigl(2a-(n-1)c\pm\sqrt{\Delta_n}\Bigr)/2,\\
  \lambda_{n;k} &:=& c+2a\,\cos(k\varphi_n),\quad k=1,2,....,n-1.\label{lamnk}
\end{eqnarray}
\end{subequations}
For higher transparency, we omit the index~$n$, so $M\equiv M_n$, 
$\beta_{n;{-}}\equiv\beta_-$ etc. With above denotations there holds true
\begin{thm}\label{egpair}
For arbitrary $n>2$ the eigenpairs of the matrix~$M$ are
\begin{equation}\label{mnep} 
      (\lambda_-,\bm{w}_-),\, (\lambda_+,\bm{w}_+)\quad\mathrm{and}\quad
     (\lambda_k,\bm{w}_k),\; k=1,2,...,n-1,
\end{equation}
where
\begin{subequations}\label{mnev}
\begin{eqnarray}
\bm{w}_\pm &=& [\beta_\pm, 1, 1, 1, \ldots, 1, 1]^\mathsf{T},
  \label{mnevpm}\\
\bm{w}_k &=& 
  [0,1,\omega^k,\omega^{2k},\ldots,\omega^{(n-2)k},\omega^{(n-1)k}]^\mathsf{T}, 
  \label{mnevk}
\end{eqnarray}
\end{subequations}
i.e., $(\bm{w}_k)_0=0$ and $(\bm{v}_k)_j=\omega^{(j-1)k}$, for 
$j=1,2,\ldots,n$.
\end{thm}

\begin{prf}
The proof consists in demonstrating that there hold true the equalities
$M\bm{w}_\pm=\lambda_\pm \bm{w}_\pm$ and $M\bm{w}_k=\lambda_k\bm{w}_k$, and 
that the collection 
$\{\bm{w}_-,\bm{w}_+,\bm{w}_1,\bm{w}_2,\ldots,\bm{w}_{n-1}\}$ is linearly 
independent. First we state that there exists a~number~$\beta$ (and we specify 
it later) such that 
\[
  \bm{w} := [\beta, 1, 1,\ldots, 1]^\mathsf{T} 
\]
is an eigenvector of~$M$. For every~$\beta$ there is 
\[
   M\bm{w}=\bigl[ nb-nc\beta,\lambda,\lambda,\ldots,\lambda\bigr]^\mathsf{T},
\]
where $\lambda=b\beta+2a+c$, so~$\beta$ is the eigenvalue of~$M$ iff 
\[
  nb-nc\beta = \beta\lambda=\beta(b\beta+2a+c). 
\]
This condition is the quadratic equation in~$\beta$, with solutions $\beta_-$ 
and~$\beta_+$ (they are distinct because $b\ne0$). For them the 
eigenvalue~$\lambda$ assumes the value $\lambda_-$ and~$\lambda_+$, 
respectively, and the vector~$\bm{w}$ is $\bm{w}_-$ and~$\bm{w}_+$. Since the 
discriminant~$\Delta$ of the equation at hand is positive for $b\ne0$, so 
$\lambda_-<\lambda_+$. The eigenvectors corresponding to them have identical 
coordinates (up to a constant multiplier each coordinate is equal to~1) but 
the first one: up to the same multiplier this coordinate is $\beta_-$ 
and~$\beta_+$, respectively.

Now we go to show that~$M$ has $n-1$ eigenvectors such that their first 
coordinate is~0. We will see even more: these eigenvectors of~$M$ are of form
\[ 
   \bm{w}_k:=[0,1,\omega^k,\omega^{2k},\ldots,\omega^{(n-1)k}]^\mathsf{T}, 
\]
where $k = 1,2,\ldots,n-1$. The index~$k$ fixed, we have 
$\bm{w}_k=[0,1,z,z^2,\ldots,z^{n-1}]^\mathsf{T}$, where $z:=\omega^k$. Then, by 
the properties of the roots of~unity,
\[
M\bm{w}_k
=\left[\begin{array}{c}
 b(1+z+z^2+\cdots+z^{n-1})\\
 c+a(z+z^{n-1}) \\ \vdots \\
 cz^{j-1}+a(z^{j-2}+z^j) \\ \vdots \\
 c z^{n-1}+a(1+z^{n-2})
\end{array}\right]
=\bigl(c+a(z+z^{-1})\bigr)\left[\begin{array}{c}
  0 \\ 1 \\ \vdots \\ z^{j-1} \\ \vdots \\ z^{n-1} 
\end{array}\right]
=\lambda_k\bm{w}_k
\]
The above proves that eigenpairs of~$M$ are $(\lambda_-,\bm{w}_-)$, 
$(\lambda_+,\bm{w}_+)$ and $(\lambda_k,\bm{w}_k)$, where $k=1,2,\ldots,n-1$. 

Eigenvectors $\bm{w}_k$, $k=1,2,\ldots,n-1$, are linearly independent. Really, 
by neglecting their first coordinate (recall, it is~0) we turn them into 
vectors $\bm{v}_k:=[1,\omega^k,\omega^{2k},\ldots,\omega^{(n-1)k}]^\mathsf{T}$, 
$k=1,2,\ldots,n-1$. They and the vector $\bm{v}_0:=[1,1,1,\ldots,1]^\mathsf{T}$ 
form the set $\{\bm{v}_0,\bm{v}_1,\bm{v}_2,\ldots,\bm{v}_{n-1}\}$. The matrix, 
whose columns are these $n$~vectors, is the Vandermonde matrix 
$V_n=[\bm{v}_0|\bm{v}_1|\bm{v}_2|\ldots|\bm{v}_{n-1}]
=V_n(1,\omega,\omega^2,\ldots,\omega^{n-1})$. Since~$V_n$ is nonsingular, its 
columns are linearly independent.

Recalling that~$\bm{v}_0$ is the reduced~$\bm{w}_+$ or~$\bm{w}_-$ (obtained by 
throwing away its first coordinate~$\beta_\pm$), and that~$\bm{w}_-$ 
and~$\bm{w}_+$ are linearly independent, we conclude that vectors 
$\bm{w}_-,\bm{w}_+,\bm{w}_1,\bm{w}_2,\ldots,\bm{w}_{n-1}$ are linearly 
independent. \qed
\end{prf}

\begin{figure}
\begin{center}
\includegraphics[width=13cm]{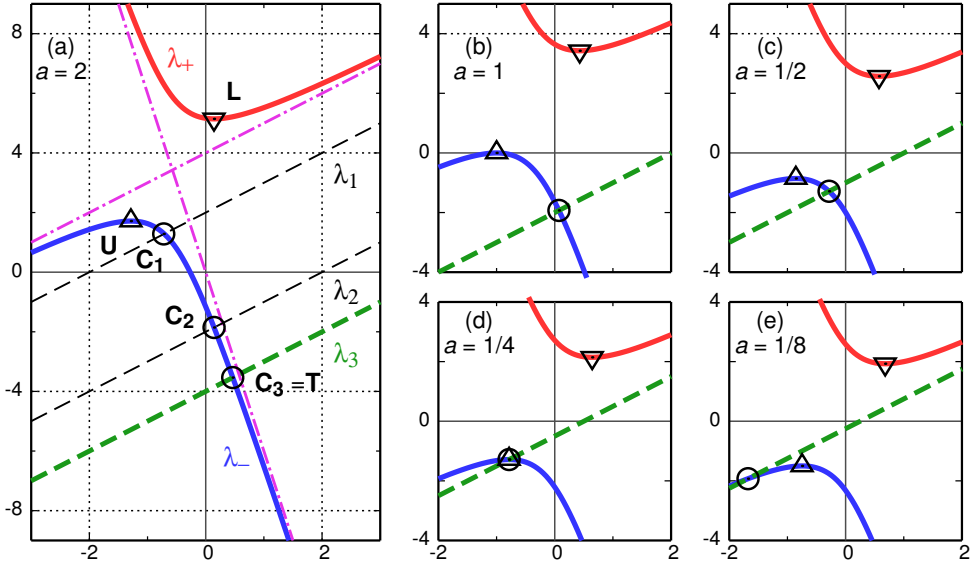}
\end{center}
\caption{(Color online)\label{m6plot}
Eigenlines $\lambda=\lambda(c)$ of an abc matrix $m_6(a,1,c)$ for (a)~$a=2$, 
(b)~$a=1$, (c)~$a=1/2$, (d)~$a=1/4$, and (e)~$a=1/8$. In Fig.~\ref{m6plot}(a) 
there are drawn seven curves (lines $\lambda_1=\lambda_5=c+a$ and 
$\lambda_2=\lambda_4=c-a$ are double), the dash-dotted lines are asymptotes to 
$\lambda=\lambda_\pm(c)$, namely $\lambda=c+2a$ and $\lambda=-nc$. There are 
also marked the uppermost (\textsf{U}) and the lowermost (\textsf{L}) points. 
Empty circles denote points~$\mathsf{C}_k$, where the lines $\lambda=\lambda_k$ 
cross the curve~$\lambda=\lambda_-$ with the transition 
point~$\mathsf{T}=\mathsf{C}_{6/2}$. These points run the curve 
$\lambda=\lambda_+(c)$ if $a<0$; see Section~\ref{SpecPt} for details.}
\end{figure}

Note that in produced formulas~$b$ appears only squared, so the sign of~$b$ 
has no importance. Plots in Fig.~\ref{m6plot} show how eigenvalues of~$M_n$ 
change in~$c$, when $b=1$ and the parameters~$a$ and~$n$ are fixed.

\section{Analyzing the eigenvalues of an abc matrix}\label{analyz}
In the whole section $b=\pm1$ is assumed, so in all cases $|b|=1$ and $b^2=1$. 
The results can be applied to the matrix~$M_2$ introduced by 
formula~(\ref{m12notilde}).

\subsection{The spectrum}\label{Spec}
Direct examination proves
\begin{lmt}\label{asmp}
  For $n>2$ and any fixed real numbers $a$, $b\ne 0$ 
\begin{enumerate}
  \item The functions $\lambda=\lambda_+(c)$ and $\lambda=\lambda_-(c)$ are 
    strictly convex and concave, respectively.\label{Lita}
\item The eigenvalues~$\lambda_+(c)$ and $\lambda_-(c)$ are separated by 
    $\lambda_\mathrm{sep}:=c+2a$, and $\lambda=\lambda_\mathrm{sep}$ is an 
    asymptote to~$\lambda=\lambda_\pm(c)$, when $c\to{\pm\infty}$, 
    respectively; moreover, these three curves, 
    $\lambda=\lambda_\mathrm{sep}(c)$ and $\lambda=\lambda_\pm(c)$, have no 
    common points.\label{Litb}
 \item The other asymptote to the curves $\lambda=\lambda_\pm(c)$ is the line 
    $\lambda=-nc$, independent of~$a$.\label{Litx}
\end{enumerate}
\end{lmt}

Since $m_n(0,b,c)=A_n(-nc,b,c)$, so, by Theorem~\ref{ahme} and 
Corollary~\ref{specA}, we have 
\begin{cor}\label{spec2}
  For $n>0$ and any real numbers~$c$ and $b\ne 0$ 
\begin{equation}\label{spec2e}
\left|\sigma\bigl(m_n(0,b,c)\bigr)\right|= \min\{n+1,3\}
\end{equation}
and for $n>1$ there is $\lambda_k=c=\lambda_\mathrm{sep}(c)$, 
$k=1,2\ldots,n-1$. 
\end{cor}

We start consideration of a~general case, $a\ne 0$, with the following
\begin{lmt}\label{specsep}
  For fixed $n>2$, a~real number $a\ne0$ and $b=\pm1$ there holds true:
\begin{enumerate}
  \item For every $k=1,2,\ldots,n-1$ the line $\lambda=\lambda_\mathrm{sep}(c)$ 
    lies above (below) $\lambda_k(c)$ for $a>0$ ($a<0$, 
    respectively).\label{Litc}
  \item There are $p_n\equiv p:=\lfloor(n-1)/2\rfloor$ pairs of equal numbers:
    $\lambda_k=\lambda_{n-k}$, for $k=1,2,\ldots,p$.\footnote{For 
    $x\in\mathbb{R}$ we use the conventional symbol 
    \protect{$\lfloor x\rfloor$} to 
    denote the greatest integer less than or equal to~$x$.}\label{Litd}
  \item The sequence $(\lambda_1,\lambda_2,\ldots,\lambda_q)$, 
    $q_n\equiv q:=\lfloor n/2\rfloor$, is strictly decreasing (increasing) for 
    $a>0$ ($a<0$, respectively). They are bounded by 
    $\lambda_\mathrm{lim}:=c-2a$, which is equal to $\lambda_q=\lambda_{n/2}$ 
    for an even number~$n$.\label{Lite}
 \item For $a>0$ ($a<0$) each $\lambda=\lambda_k(c)$ crosses 
    $\lambda=\lambda_-(c)$ ($\lambda_+(c)$, respectively) at the unique point 
    $\mathsf{C}_k$ with the abscissa 
\begin{equation}\label{ck}
  c_k=
\frac{4a^2\cos(k\varphi)(1-\cos(k\varphi))+n}{2(n+1)a(\cos(k\varphi)-1)}.
\end{equation}
    When~$k$ increases from~1 to~$q$, these abscissas form an increasing 
    (a~decreasing) sequence for $a>0$ ($a<0$, respectively).\label{Litf}
\end{enumerate}
\end{lmt}
\begin{prf}
These four claims follow properties of the cosine function. For 
$1\le k\le q$ the angles $k\varphi$ form the increasing sequence and 
$0<k\varphi\le\pi$, so the sequence $\bigl(\cos(k\phi)\bigr)_{k=1,2,\ldots,q}$ 
is decreasing, and this confirms~(\ref{Lite}). The point~(\ref{Litc}) is 
valid due to the upper bound $\cos(k\varphi)<1$ for $k\varphi>0$, whereas 
(b)~is proved due to the parity of the cosine function, 
$\cos(k\varphi)=\cos\bigl((n-k)\varphi\bigr)$. The numbers~$c_k$ are unique 
solutions of the appropriate equations
\[
 c + 2a \cos(k\varphi) = \left(2a-(n-1)c\pm\sqrt{\Delta}\right)/2.
\]
Due to the monotonic behavior of the sequence 
$\bigl(\lambda_k\bigr)_{k=1,2,\ldots,q}$, the sequence 
$\bigl(c_k\bigr)_{k=1,2,\ldots,q}$ is also monotonic and its character follows 
the properties of the curves $\lambda=\lambda_\pm$ and the sign of the 
parameter~$a$. This completes the proof of~(\ref{Litf}), so of the whole lemma, 
as well.
\qed
\end{prf}

\begin{thm}[Cardinality of the spectrum]\label{cardsig} With the same 
assumptions as in Lemma~\ref{specsep} we have
\[
  |\sigma(M)|=\begin{cases}
        \lfloor n/2\rfloor+1, & 
             \mathrm{if}\qquad c=\in\{c_1,c_2,\ldots,c_q\},\\
        \lfloor n/2\rfloor+2, & \mathrm{otherwise}.
  \end{cases}
\]
\end{thm}
\begin{prf} This theorem follows immediately the claims in Lemma~\ref{specsep}.
\qed
\end{prf}
Note that for an even number~$n$ there is $q=p+1=n/2$ and, for every 
$c\notin\{c_1,c_2,\ldots,c_q\}$, the spectrum of~$M$ has $n/2+2$ elements, with 
three single eigenvalues ($\lambda_\pm$ and $\lambda_{n/2}$) and $n/2-1$ double 
ones. When the number $c=c_k$, with $k<n/2$, there are two single eigenvalues, 
$(n/2-2)$ double eigenvalues, and one triple eigenvalue. With $c=c_{n/2}$ these 
numbers are~1, $n/2$, and~0, respectively. If the number~$n$ is odd, then 
$q=p=(n-1)/2$ and there is no single eigenvalue~$\lambda_k$. Therefore, there 
are $(n+3)/2$ eigenvalues in a~general case (two single and $(n-1)/2$ double 
ones) and $(n+1)/2$ of them for $c=c_k$ (there are~1, $(n-3)/2$, and~1 
eigenvalues of the multiplicity~1, 2, and~3, respectively).

\subsection{The special points}\label{SpecPt}
Since the function $\lambda=\lambda_-(c)$ is strictly concave in~$c$, then it 
may have the (local, so also global) maximum and, in fact, it has. Its 
uppermost point is $\mathsf{U}_n(a):=(c_\mathrm{upp},\lambda_\mathrm{upp})$, 
where
\begin{subequations}\label{upper}
\begin{eqnarray}
  c_\mathrm{upp}(a,n) := -\frac{(n-1)+2a}{n+1} = -1-\frac{2(a-1)}{n+1}, 
   \label{uppc}\\
   \lambda_\mathrm{upp}(a,n)  :=\frac{2n}{n+1}(a-1)=-n(c_\mathrm{upp}+1).
\label{uppl}
\end{eqnarray}
\end{subequations}
Similarly, the convex function $\lambda=\lambda_+(c)$ has the lowermost point 
(the global minimum) at 
$\mathsf{L}_n(a):=(c_\mathrm{low},\lambda_\mathrm{low})$, with
\begin{subequations}\label{lower}
\begin{eqnarray}
  c_\mathrm{low}(a,n)  := \frac{(n-1)-2a}{n+1} = 1-\frac{2(a+1)}{n+1}, 
   \label{lowc}\\
   \lambda_\mathrm{low}(a,n)  :=\frac{2n}{n+1}(a+1)=-n(c_\mathrm{low}-1).
 \label{lowl}
\end{eqnarray}
\end{subequations}

Considerations of the other ``special point'' we start with 
\begin{cor}
For fixed $a\ne0$ and any number $n>2$ the abscissas $c_k$, 
$k=1,2,\ldots,q$, belong to the interval $(-\infty,c_-]$, for $a>0$, and to 
the interval $[c_+,\infty)$, for $a<0$, where $c_\pm$ denote the abscissa of 
the point at which $\lambda_\pm(c)=\lambda_\mathrm{lim}(c)$.
\end{cor}
\begin{prf}
At first we note that for arbitrary fixed $a\ne0$, $n>2$, and $k=1$ 
Eq.~(\ref{ck}) gives
\[
   c_1=
\frac{-2a\cos\varphi}{n+1} 
   + \frac{n}{n+1}\, \frac{1}{2a(\cos\varphi-1)}.
\]
Obviously, the first summand tends to~0. Moreover, $\cos\varphi<1$ and 
$\lim_{n\to\infty}(\cos\varphi-1)=0$, therefore $\lim_{n\to\infty}c_1=-\infty$ 
for $a>0$ and $\lim_{n\to\infty}c_1=\infty$, otherwise. The existence and 
properties of the points~$c_\pm$ follow from the points~(\ref{Lite}) 
and~(\ref{Litf}) in Lemma~\ref{specsep}, so the proof is completed.\qed
\end{prf}

In the other domains, i.e.\ for $c\in(c_-,\infty)$, for $a>0$, and 
$c\in(-\infty,c_+)$, for $a<0$, the graphs $\lambda=\lambda_\pm(c)$ and 
$\lambda=\lambda_k(c)$ have no common points, so in these domains 
$|\sigma(M)|=q+2$. The spectrum cardinality at~$c_\pm$ depends on the parity 
of~$n$, namely $|\sigma(M)|=\lfloor(n+3)/2\rfloor$.

Taking into account Eq.~(\ref{ck}) we see that both numbers $c_\pm$ are 
determined by the same formula. We refer to the point, at which the one of the 
curves~$\lambda=\lambda_\pm(c)$ meets the limit line 
$\lambda=\lambda_\mathrm{lim}$, as a~\textit{transition point\/} (generated by 
$a\ne0$ and~$n$). Coordinates of this point, 
$\mathsf{T}_n(a):=(c_\mathrm{trans},\lambda_\mathrm{trans})$, can be determined 
substituting $\cos(k\varphi)=-1$ in Eq.~(\ref{ck}) and then taking 
$\lambda_\mathrm{trans}=c_\mathrm{trans}-2a$. In this way we have
\begin{subequations}\label{trans}
\begin{eqnarray}
  c_\mathrm{trans}(a,n) &=& \frac{8a^2-n}{4(n+1)a},\label{transc}\\
  \lambda_\mathrm{trans}(a,n) &=& -\frac{n(8a^2+1)}{4(n+1)a}.\label{transl}
\end{eqnarray}
\end{subequations}
The transition point $\mathsf{T}_n(a)$ sits on the curve 
$\lambda=\lambda_-(c)$, for $a>0$, and on the curve $\lambda_+(c)$ for $a<0$.

\begin{rmk}\label{trans2}
Note that Eqs.~(\ref{trans}) for $n=2$ give
\begin{equation}\label{T2notil}
  \mathsf{T}_2(a)=\frac{1}{6 a}\left(4a^2-1,-(8a^2+1)\right)
\end{equation}
and this formulas correspond to the matrix~$M_2$ given by 
formula~(\ref{m12notilde}). To obtain them for the matrix~$\widetilde{M}_2$, 
determined in by~(\ref{m12tilde}), we have to replace each $2a$ by~$a$, so
(cf.~Theorem~\ref{thmtilde}) 
\begin{equation}\label{T2til}
  \widetilde{\mathsf{T}}_2(a)=\frac{1}{3 a}\left(a^2-1,-(2a^2+1)\right).
\end{equation}
\end{rmk}.

\begin{figure}
\begin{center}
\includegraphics[width=10cm]{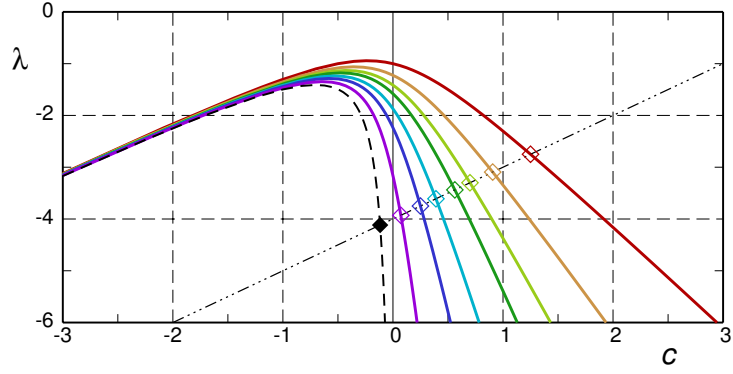}
\end{center}
\caption{(Color online)\label{trlim}
Transition curves $\mathcal{T}_n(a)$, for $a>0$, i.e., curves 
$(c,\lambda) = T_n(a)$, traced for $b=1$ and $n=2,3,4,5,7,10,20$ (solid 
lines from right to left), the limit transition curve $\mathcal{T}(a)$ 
(dashed line) and the line $\lambda_\mathrm{lim}=c-2a$ for $a=2$ (dash-dotted 
line) are also plotted. All the transition points $T_n(2)$ (empty diamonds) sit 
on this line. This line also keeps the limit point 
$T_\infty(2)=\lim_{n\to\infty}T_n(2)=(-1,-33)/8$ (full diamond), it is the 
point where the curve $\mathcal{T}(a)$ and the line $\lambda=c-2a$ intersect, 
for fixed $a=2$. In the case $n=2$ the matrix~$M_2$ is taken into account.}
\end{figure}

For fixed numbers~$n$ the set $\{T_n(a)\mid a\ne 0\}\equiv
   \{(c_\mathrm{trans},\lambda_\mathrm{trans})\}\subset\mathbb{R}^2$ determines 
the graph~$\mathcal{T}_n(a)$ hereafter referred to as the $n$th
\textit{transition curve}. Note that this curve has two separate branches: one 
for $a>0$ and the second for $a<0$. These branches have asymptotes: 
$\lambda=c$ and $\lambda=-nc$. Sending~$n$ to infinity we get
\begin{equation}\label{limtp} 
       T_\infty(a) := \lim_{n\to\infty} T_n(a) = \frac{-1}{4a}(1,8a^2+1). 
\end{equation}
This point referred to as a~\textit{limit transition point\/} (associated 
to~$a$); the set (and the graph) $\mathcal{T}(a)=\{T_\infty(a)\mid a\ne 0\}$ is 
called a~\textit{limit transition curve}. This limit curve is explicitly 
described via the relation $\lambda= c+1/(2c)$, $c\ne 0$. For $c<0$ ($c>0$), so 
also $\lambda<0$ ($\lambda>0$, respectively), this graph is the ``limit curve'' 
for the lower (the upper) branch of $\mathcal{T}_n(a)$, i.e., for $a>0$ and 
$a<0$, respectively. Since all the transition points satisfy the relation
$\lambda=\lambda_\mathrm{lim}(c)=c-2a$ and, moreover, it is also satisfied by 
points sitting on the graph~$\mathcal{T}(a)$, then at each point~$T_\infty(a)$ 
the curve $\mathcal{T}(a)$ and the line $\lambda_\mathrm{lim}$ intersect each 
other. Fig.~\ref{trlim} shows some transition curves and the limit transition 
curve.

\subsection{The extreme eigenvalues}\label{Extr}
The points $T_n(a)$ are actual transitions points for an even number~$n$, so,  
hereafter, we restrict ourselves to this case. Solving equations (for $a>0$)
\[
 c_\mathrm{upp}=c_\mathrm{trans}\qquad\mathrm{and}\qquad
  \lambda_\mathrm{upp}=\lambda_\mathrm{trans}
\]
we find the value $a_{\mathrm{crit}+}=1/4$ at which the transition point 
$\mathsf{T}_n(a)$ coincides with the uppermost point $\mathsf{U}_n(a)$. The 
same procedure gives $a_{\mathrm{crit}-}=-1/4$ for which 
$\mathsf{T}(a)=\mathsf{L}(a)$ ($a<0$). Since $\lambda_\mathrm{upp}$ 
($\lambda_\mathrm{low}$) is the global maximum (minimum, respectively), then 
for $a\ne a_{\mathrm{crit}\pm}$ we always have 
$\lambda_\mathrm{trans}<\lambda_\mathrm{upp}$, when $a>0$, and 
$\lambda_\mathrm{trans}>\lambda_\mathrm{low}$, when $a<0$. Now, directly from 
Lemmas~\ref{asmp} and~\ref{specsep} it follows
\begin{cor}\label{abobel} The following inequalities are satisfied (for $n>2$)
\begin{eqnarray*} 
 c_\mathrm{trans} &<& c_\mathrm{low},\qquad \mathrm{for}\quad a<-1/4,\\
 c_\mathrm{trans} &>& c_\mathrm{low},\qquad \mathrm{for}\quad -1/4<a<0,\\
.c_\mathrm{trans} &<& c_\mathrm{upp},\qquad \mathrm{for}\quad 0<a<1/4,\\
 c_\mathrm{trans} &>& c_\mathrm{upp},\qquad \mathrm{for}\quad a>1/4.
 \end{eqnarray*}
The critical values $a_{\mathrm{crit}\pm}=\pm1/4$ do not depend on~$n$.
\end{cor}
These relations are also satisfied for the matrix~$M_2$. When the 
matrix~$\widetilde{M}_2$ is at hand, then 
$\tilde{a}_{\mathrm{crit}\pm}=2a_{\mathrm{crit}\pm}=\pm1/2$. 

Reassuming the above, we have
\begin{thm}[Extreme eigenvalues]\label{mainthm}
For any natural number $n>2$ the extreme eigenvalues $\lambda_\mathrm{min}$ and 
$\lambda_\mathrm{max}$ of the matrix $m_n(a,b=\pm1,c)$ are given by the 
following formulas
\begin{description}
 \item [Case $\bm{a<0}$:]
\begin{equation}\label{mainneg}
 \lambda_\mathrm{min}(c)=\lambda_-(c),\qquad\mathrm{and}\qquad
 \lambda_\mathrm{max}(c)=\begin{cases} 
    \lambda_+(c), & \mathrm{for}\quad c< c_\mathrm{trans},\\
    \lambda_{n/2}=c-2a,  & \mathrm{for}\quad c\ge c_\mathrm{trans};
 \end{cases}
\end{equation}
 \item [Case $\bm{a=0}$:] $\lambda_\mathrm{min}(c)=\lambda_-(c)$ and 
   $\lambda_\mathrm{max}(c)=\lambda_+(c)$;
  \item [Case $\bm{a>0}$:] 
\begin{equation}\label{mainpos}
 \lambda_\mathrm{min}(c)= \begin{cases} 
  \lambda_{n/2}=c-2a,  & \mathrm{for}\quad c\le c_\mathrm{trans},\\
 \lambda_-(c) & \mathrm{for}\quad c> c_\mathrm{trans},
 \end{cases}
\qquad\mathrm{and}\qquad
\lambda_\mathrm{max}(c)=\lambda_+(c).
\end{equation}
\end{description}
\end{thm}

The preceding claims prove the final
\begin{cor}[Extreme extrema]\label{final}
The functions $\lambda=\lambda_\mathrm{min}(c)$ and 
$\lambda=\lambda_\mathrm{max}(c)$ have their extrema 
$\min_{c\in\mathbb{R}}\lambda_\mathrm{max}(c)$ and 
$\max_{c\in\mathbb{R}}\lambda_\mathrm{min}(c)$ at
\begin{subequations}\label{fineq}
\begin{eqnarray}
 (c_\mathrm{min-of-max}, \lambda_\mathrm{min-of-max}) &=& \begin{cases}
(c_\mathrm{trans},\lambda_\mathrm{trans}), 
     & \mathrm{for} \quad a\le -1/4\\
(c_\mathrm{low},\lambda_\mathrm{low}), 
     & \mathrm{for} \quad a>-1/4;
\end{cases} \label{minmax}\\
(c_\mathrm{max-of-min}, \lambda_\mathrm{max-of-min}) &=& \begin{cases}
(c_\mathrm{trans},\lambda_\mathrm{trans}), 
     & \mathrm{for} \quad a\ge 1/4\\
(c_\mathrm{upp},\lambda_\mathrm{upp}), 
     & \mathrm{for} \quad a<1/4.
\end{cases} \label{maxmin} \label{phys}
\end{eqnarray}
\end{subequations}
Whenever such extremum coincides with the transition point, the corresponding 
eigenvalue (the minimum of $\lambda_\mathrm{max}$ or the maximum of 
$\lambda_\mathrm{min}$) is double, and it is single, otherwise.
\end{cor}

\section{An abc matrix and its corresponding wheel graph}\label{wheel}
\begin{figure}
\begin{center}
\includegraphics[width=12cm]{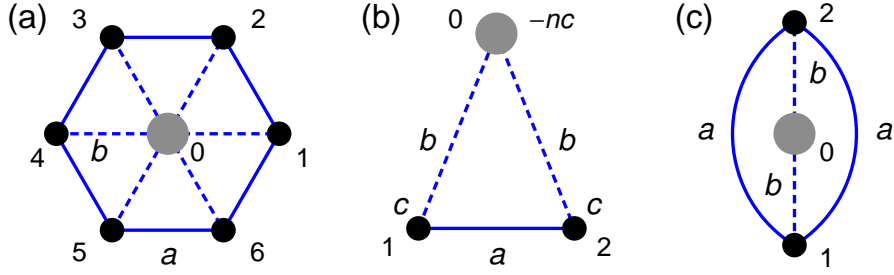}
\end{center}
\caption{(Color online)\label{wngraphs}
(a) The wheel graph $W_7\equiv W_{6+1}$ (the hub or the central vertex is 
labeled by~0, and tire vertices are labeled by $1, 2,\ldots,n$, ) with its 
edge-weights (the tire edges of weight~$a$, solid lines, and the spokes of 
weight~$b$, dashed lines). For $n=2$ (b) a~weighted cyclic graph~$C_3$, with 
the adjacency matrix~$\widetilde{M}_2$, or (c) a~multigraph, with the adjacency 
matrix~${M}_2$, can be considered. In all cases black full circles denote 
vertices of weight~$c$, whereas gray (and larger) ones denote those of 
weight~$-nc$. 
}
\end{figure}

Just introduced abc matrix~$m_n(a,b,c)$ can be interpreted as the weighted 
adjacency matrix of a~regularly weighted wheel graph, the notion we define 
below and we mostly follow the nomenclature presented by \citet{Bond76} and 
\citet{Bran99} \citep[see also][]{WolfA}. To unify the nomenclature, we say 
that an $N$-vertex graph, $N>3$, in which its $n:=N-1$ vertices (referred to as 
\textit{cyclic}, or \textit{tire vertices}) form the cyclic graph~$C_n$ and one 
single vertex is universal (i.e., is adjacent to every vertex), is called 
an $N$th \textit{wheel (graph)}, and denoted by~$W_N$.\footnote{Some authors 
\citep[e.g.,][]{Rose11} use $W_n$ to denote this graph.\label{Nnidx}} One can 
notice that $W_N=K_1+C_n$,where $K_1$~is the singleton graph \citep{Ranj10}. 
The universal vertex is also called a~\textit{hub}, or a~\textit{central 
vertex}, and every edge adjacent to it is called a~\textit{spoke}. $W_N$ is 
a~planar graph and can be seen as a~centered (regular) $n$-gon or as 
a~skeleton of the (regular) pyramid whose base is an $n$-gon, then its 
universal vertex is simply its apex, tire edges and spokes are base and lateral 
edges, respectively; in this context one can call $W_N$ an $n$th pyramid graph. 
An example, the $W_7$ graph, is presented in Fig.~\ref{wngraphs}(a).

When to each of $2n$~edges and to each vertex $j=0,1,\ldots,n$ of~$W_N$ there 
are assigned some real numbers (referred to as weights), we have 
a~\textit{weighted wheel graph}, also denoted by~$W_N$. It is said to be 
\textit{regular(ly)}, if: (a)~every edge $(j,j+1)$, $n+1\equiv 1$, in~$C_n$ is 
of weight~$w(j,j+1):=a$; (b)~every spoke $(0,j)$ is of weight~$w(0,j):=b$; 
(c)~every tire vertex is of the same weight~$w(j)=c$; (d)~the weight of the 
central vertex is~$w(0)=-nc$. Here~$b$ can be seen as a~scale, and without the 
loss of generality we can take into account the sign of~$b$ only, i.e., we 
restrict discussion to the cases $b=\pm1$ (cf.\ Remark~\ref{trabc}). Clearly, 
$W_N$~is of dihedral symmetry~$D_n$, although this group does not act 
transitively on all vertices. In a~more general case, it can be assumed 
$w(0)=c'\ne-nc$, what preserves the dihedral symmetry. However, in physical 
problems the matrices~$M_n$ are applied to, one needs traceless matrices (see 
further in the text for details), so we demand $w(0)=-nc$ in a~regular weighted 
wheel graph~$W_N$. Such graphs, and corresponding weighted adjacency matrices, 
are considered in some papers representing various fields of science 
\citep[see, e.g.,][]{Bapa10,Stev11,Patt14,CaiH17}. Note that for $a=0$ we have 
no wheel graph~$W_N$ any longer (since the edges with zero weights are 
removed), but we have a~star graph~$S_N$ (text of footnote~\ref{Nnidx} applies 
accordingly), which is connected. This is why we have chosen the weight of 
spokes not of the tire (cyclic) edges to be a~scale in problems discussed here.

In the main text there are considered matrices~$\widetilde{M}_2$ and $M_2$. 
Formally, they \textit{are not\/} weighted adjacency matrices of wheel graphs, 
since usually it is assumed that a~wheel graph~$W_N$ has at least four 
vertices. By formula~(\ref{m12tilde}) it is clear that the 
matrix~$\widetilde{M}_2$ is a~weighted adjacency matrix of the cyclic 
graph~$C_3$ with one edge of different weight~$a$, whereas the other two edges 
are of equal weights~$b$ [see Fig.~\ref{wngraphs}(b)]. So this graph (and the 
corresponding matrix) has the $D_2$ (the Klein group) symmetry. A~doubled 
weight~$2a$ in the matrix~$M_2$ appears when we sketch the cyclic graph~$C_2$, 
with two vertices and \textit{two\/} edges, i.e., a~digon (two-gon); due to two 
undirected edges (1,2) it is a~multigraph [see Fig.~\ref{wngraphs}(c)]. 
However, it has, again, the $D_2$~symmetry and connecting its vertices to the 
hub (labeled by~0) we obtain ``a~centered digon''. The tire edges have the same 
weight~$a$, so (at least in some applications) this graph (its adjacency 
matrix, in fact) has the same properties as the cyclic graph~$C_3$ with one 
distinguished edge of weight~$2a$. This shows that the matrix~$M_2$ can be 
included in the discussion presented, whereas the matrix~$\widetilde{M}_2$ 
constitutes the special case.

\section{Final remarks}\label{cncl}
In this paper we introduced arrow-bordered circulant (``abc'' for short) 
matrices $m_n(a,b,c)$ for $n\ge2$ and real parameters $a$, $b$, and~$c$. For 
$a=0$ and $bc\ne0$ the arrowhead matrices are revealed, whereas assuming $b=0$ 
and $a\ne0$ one obtains two blocks: a~trivial one-dimensional matrix and 
a~circulant matrix with at most three nonzero elements in each row. Such 
object can be considered as weighted adjacency matrices of wheel graphs~$W_N$ 
(or star graphs~$S_N$ for $b=0$). In this paper the special case has been 
investigated: regular abc matrices (see Section~\ref{abcm}). They inherit some 
properties of their ``parents'', arrowhead and circulant matrices. We have 
determined eigenpairs of the abc matrices and discussed their eigenvalues. We 
have put stress on their bounds, asymptotic behavior, and extrema. Since such 
matrices are widely used in different fields of science, then this work 
provides results of some interest. It is desirable to have results in more 
general cases, e.g., for irregular abc matrices, when the vector~$\bm{b}$ is 
not constant \citep[cf.,][]{Shen09}, or arrow-bordered alternating circulant 
matrices \citep{Tee07}. In the latter case signs of consecutive rows are 
alternated, but in more general cases rows with alternated \textit{values\/} of 
nonzero matrix elements should also be included.

This work has been motivated by some physical problems, in particular 
these related to magnetic (finite) systems. Classical counterpart of the 
Heisenberg model describes a~set of \textit{localized\/} spin 
vectors~$\bm{s}_j$ with interactions determined by real numbers~$J_{ij}$, where 
``locations'' of vectors~$\bm{s}_j$ are labeled by $i,j=1,2,\ldots,N$; this 
numbers are interpreted as graph vertices (to each vertex~$j$ of a~given 
graph~$G$ a~spin vector~$\bm{s}_j$ is assigned). Within this approach couplings 
$J_{ij}$ are nondiagonal elements of a~weighted adjacency matrix (with 
vanishing diagonal) of the graph~$G$. The total energy of this system equals 
$E=\sum_{(j,j')}J_{jj'}\bm{s}_j\cdot\bm{s}_{j'}$, where the standard inner 
product is denoted by dot~``$\cdot$''. \citet{Schm03} \citep[see also a~series 
of papers][]{Schm1,Schm2,Schm3,Schm4} introduced the so-called gauge vector 
$\bm{c}\in\mathbb{R}^N$ and then they ``dressed'' the matrix~$J$ assuming 
$J_{jj}=c_j$, so each number~$c_j$ can be considered as weight of the 
vertex~$j$. They showed that for traceless dressed matrix (i.e., for 
$\sum_j c_j=0$) some physical quantities are ``gauge independent''. For fixed 
parameters of the system the eigenvalues of the dressed matrix, including 
its minimum eigenvalue $\lambda_\mathrm{min}(\bm{c})$, depend on the gauge 
vector~$\bm{c}$. They proved that there exists 
$\max_{\bm{c}\in\mathbb{R}^N}\lambda_\mathrm{min}(\bm{c})$ and it equals, up to 
a~constant factor, the minimum of the system energy~$E$. The considerations 
performed in this paper give us a~general solution of this problem for spins 
placed at vertices and the center of a~regular $n$-gon with an even 
number~$n$ \citep[see, e.g.,][]{Graj18,FAK,Ako07,Kaka18}. The further 
physical analysis take into account also eigenvectors of the dressed 
matrix~$J$, but this problem is out of scope here. It has to be emphasized that 
the degeneracy of the determined maximum of the minimum eigenvalue 
$\lambda_\mathrm{min}(\bm{c})$ says whether the lowest energy configuration of 
spin vectors is collinear, coplanar or spatial---it happens 
for single, double and triple~$\lambda_\mathrm{min}(\bm{c})$ (cf.\ remarks at 
the end of Section~\ref{Spec}). Note that higher degeneracy is not excluded and 
\citet{Schm1} have provided an example of a~classical spin system, when the 
lowest energy configuration can be realized in nonphysical four-dimensional 
space~$\mathbb{R}^4$. For small systems some results have been obtained with 
simple calculus \citep{Graj18,KFA}, but more general considerations need 
strictly proved properties, i.e., the results of this paper \citep[see][]{FMK}. 
It is worth noting that for actually synthesized magnetic molecules 
\citep[see, e.g.,][]{Bani18, Maje18, Prsa18} there is needed analysis of some 
more general matrices (not regularly weighted graphs), e.g., alternating circulant 
matrices or arrow-bordered alternating circulant ones.



\bibliography{abcm}

\end{document}